\documentclass[11pt]{article}
\usepackage{amsmath}
\usepackage{amssymb}
\usepackage{amsthm}
\usepackage{url}

\numberwithin{equation}{section}
\newtheorem{theorem}{Theorem}[section]
\newtheorem{proposition}[theorem]{Proposition}
\newtheorem{Remark}[theorem]{Remark}
\newenvironment{remark}{\begin{Remark}\rm}{\end{Remark}}
\newcommand\Proof{\noindent{\bf Proof}\quad}

\newcommand\nonu{\nonumber}
\newcommand\dstyle\displaystyle
\newcommand\sa{\smallskipamount}
\newcommand\ma{\medskipamount}
\newcommand\ba{\bigskipamount}
\newcommand\sLP{\\[\sa]}
\newcommand\mLP{\\[\ma]}
\newcommand\bLP{\\[\ba]}
\newcommand\bPP{\\[\ba]\indent}
\newcommand\RR{\mathbb{R}}
\newcommand\ZZ{\mathbb{Z}}
\newcommand\FSE{{\cal E}}
\newcommand\FSJ{{\cal J}}
\newcommand\goH{\mathfrak{H}}
\newcommand\al\alpha
\newcommand\be\beta
\newcommand\de\delta
\newcommand\la\lambda
\newcommand\tha\theta
\newcommand\Ga{\Gamma}
\newcommand\half{\frac12}
\newcommand\thalf{\tfrac12}
\newcommand\iy\infty
\newcommand\lan{\langle}
\newcommand\ran{\rangle}
\newcommand\wh{\widehat}
\newcommand\wt{\widetilde}

\newcommand\RHS{right-hand side}
\newcommand\Equiva{\Longleftrightarrow}
\newcommand{\hyp}[5]{\,\mbox{}_{#1}F_{#2}\!\left(
  \genfrac{}{}{0pt}{}{#3}{#4};#5\right)}
\newcommand{\qhyp}[5]{\,\mbox{}_{#1}\phi_{#2}\!\left(
  \genfrac{}{}{0pt}{}{#3}{#4};#5\right)}
\newcommand\bma{\begin{pmatrix}}
\newcommand\ema{\end{pmatrix}}
\newcommand\DAHA{\tilde\goH}

\begin{document}

\title{Nonsymmetric Askey-Wilson polynomials as vector-valued
polynomials}
\author{Tom H. Koornwinder\footnote{Korteweg-de Vries Institute,
University of Amsterdam, P.O.~Box 94248, 1090 GE Amsterdam, Netherlands;
\tt T.H.Koornwinder@uva.nl}\;
and Fethi Bouzeffour\footnote{Institut Pr\'eparatoire aux \'Etudes
d'Ing\'enieur de Bizerte, 8000 Bizerte, Tunisia;
\tt bouzeffourfethi@yahoo.fr}}
\date{\em{Dedicated to Paul Butzer on the occasion of his 80th birthday}}
\maketitle
\begin{abstract}
Nonsymmetric Askey-Wilson polynomials are usually written as Laurent
polynomials. We write them equivalently as 2-vector-valued symmetric
Laurent polynomials. Then the Dunkl-Cherednik operator of which they are
eigenfunctions, is represented as a $2\times 2$ matrix-valued operator.
As a new result made possible by this approach we obtain positive
definiteness of the inner product in the orthogonality relations,
under certain constraints on the parameters. A limit transition to
nonsymmetric little $q$-Jacobi polynomials also becomes possible in this way.
Nonsymmetric Jacobi polynomials are considered as limits both of
the Askey-Wilson and of the little $q$-Jacobi case.
\end{abstract}

%until 72
\section{Introduction}
Originally, most orthogonal special functions associated with root systems, like the Heckman-Opdam Jacobi polynomials and the Macdonald polynomials,
were Weyl group invariant functions. By work by Dunkl, Heckman and in
particular Cherednik it was next shown that there are related orthogonal
systems of special functions which are not Weyl group invariant, but are
in a sense more simple, and from which the earlier Weyl group invariant
special functions can be obtained by symmetrization.
While Cherednik's theory, involving double affine Hecke algebras, has been
developed for general root systems, to a large extent independent of
the classification, the specialization of this theory to the case of
rank one has its own interest, because everything can be done there
in a much more explicit way, and new results for special functions in one
variable can be obtained. In the rank one case the Weyl group
has order 2, and Weyl group symmetry turns down to a symmetry under the map
$t\mapsto-t$ or $z\mapsto z^{-1}$.

The case of nonsymmetric Askey-Wilson polynomials was first considered briefly
by Sahi \cite{6}, \cite{7} as a specialization for $l=1$ of his extension
of Cherednik's theory to the pair of root systems
$(C_l^\vee,C_l)$. The case $l=1$ was treated afterwards in more detail
by Noumi and Stokman \cite{5} and by Macdonald \cite[\S6.6]{4}.
See also the first author's paper \cite{3}.
In all these references the nonsymmetric Askey-Wilson
polynomials are defined
as certain Laurent polynomials and the bilinear form for which they
are orthogonal (or rather biorthogonal) is given in terms of a contour
integral.

In the present paper we give a presentation of nonsymmetric
Askey-Wilson polynomials as two-dimensional vector-valued polynomials.
The vector is a pair of symmetric Laurent polynomials: an Askey-Wilson
polynomial and another one with the degree lowered by 1 and two parameters
raised by 1. As a first application of this approach we can treat the
orthogonality in the nonsymmetric case by reducing it to the
well-known orthogonality of the two components of the vector. We obtain
the new result that, under certain  constraints on the parameters, there
is orthogonality with respect to a positive definite inner product. It does
not seem possible to read off this result from a positive weight function
in the contour integral.

As a second application we can easier take limits of nonsymmetric
Askewy-Wilson polynomials corresponding to the familiar limits
of the symmetric polynomials in the ($q$-)Askey scheme. In this paper
we consider the limit from nonsymmetric Askey-Wilson polynomials
to nonsymmetric little $q$-Jacobi polynomials.
Again, the vector-valued approach is crucial here.
Next we obtain nonsymmetric Jacobi polynomials as limits in two
different ways: directly as a limit of the Askey-Wilson case,
where the Laurent polynomials can be used, and as a limit of
the little $q$-Jacobi case, where we need the vector-valued polynomials.
We start the exposition by recalling the well-known nonsymmetric Bessel
case, and we end by obtaining nonsymmetric Bessel as a limit of
nonsymmetric Jacobi.

We expect that to every limit arrow in the
($q$-)Askey scheme there corresponds a limit transition in the
nonsymmetric case, and also for limits going out of the
($q$-)Askey scheme to various ($q$-)Bessel functions. We hope to present
some further examples in a subsequent paper. Another topic which we
will not yet treat here, is a study of the corresponding limit algebras of
the double affine Hecke algebra in the Askey-Wilson case.
\bLP
{\bf Conventions}\quad
Throughout assume $0<q<1$.
For \mbox{($q$-)Pochhammer} symbols and
($q$-)hyper\-geometric series
use the notation of \cite{1}. In particular, for $n\ge0$,
\begin{gather*}
(a;q)_n:=\prod_{j=0}^{n-1}(1-aq^j),\qquad
(a_1,\ldots,a_r;q)_n:=(a_1;q)_n\cdots(a_r;q)_n\,,
\\
\qhyp r{r-1}{q^{-n},a_2,\ldots,a_r}{b_1,\ldots,b_{r-1}}{q,z}:=
\sum_{k=0}^n\frac{(q^{-n},a_2,\ldots,a_r;q)_k}
{(b_1,\ldots,b_{r-1},q;q)_k}\,z^k.
\end{gather*}

Let $e_1$, $e_2$, $e_3$, $e_4$ be the elementary symmetric polynomials in $a$, $b$, $c$, $d$:
\begin{gather}
e_1:=a+b+c+d,\qquad
e_2:=ab+ac+bc+ad+bd+cd,\nonumber\\
e_3:=abc+abd+acd+bcd,\qquad
e_4:=abcd.
\label{21}
\end{gather}

For Laurent polynomials $f$
in $z$ the $z$-dependence will be written as $f[z]$.
Symmetric Laurent polynomials
$f[z]=\sum\limits_{k=-n}^n c_k z^k$ (where $c_k=c_{-k}$) are related to ordinary
polynomials $f(x)$ in $x=\thalf(z+z^{-1})$ by
$f(\thalf(z+z^{-1}))=f[z]$.
\section{The nonsymmetric Hankel transform}
\label{69}
For Bessel functions $J_\al$ see \cite[Ch. 10]{8} and references given there.
We will work with differently normalized Bessel functions
\begin{equation}
\FSJ_\al(x):=\Ga(\al+1)\,(2/x)^\al\,J_\al(x).
\label{4}
\end{equation}
Then (see \cite[(10.16.9)]{8})
\[
\FSJ_\al(x)=
\sum_{k=0}^\iy\frac{(-\tfrac14 x^2)^k}{(\al+1)_k\,k!}
=\hyp01-{\al+1}{-\tfrac14x^2}\qquad(\al>-1).
\]
$\FSJ_\al$ is an entire analytic function and we have the simple
properties and special cases
\[
\FSJ_\al(x)=\FSJ_\al(-x),\quad\FSJ_\al(0)=1,\quad
\FSJ_{-1/2}(x)=\cos x,\quad
\FSJ_{1/2}(x)=\frac{\sin x}x\,.
\]
The functions $x\mapsto\FSJ_\al(\la x)$ satisfy the eigenvalue equation
\cite[(10.13.5)]{8}:
\[
\left(\frac{d^2}{dx^2}+\frac{2\al+1}x\,\frac d{dx}\right)
\FSJ_\al(\la x)=-\la^2\,\FSJ_\al(\la x).
\]
The {\em Hankel transform} pair \cite[\S10.22(v)]{8}, for $f$ in a suitable
function class, is given by
\begin{equation}
\begin{cases}
&\dstyle\wh f(\la)=\int_0^\iy f(x)\FSJ_\al(\la x) x^{2\al+1}\,dx,\mLP
&\dstyle
f(x)=\frac1{2^{2\al+1}\Ga(\al+1)^2}\int_0^\iy \wh f(\la)
\FSJ_\al(\la x) \la^{2\al+1}\,d\la.
\end{cases}
\label{1}
\end{equation}

Now consider the so-called
{\em nonsymmetric Bessel function},
also called {\em Dunkl-type Bessel function}, in the rank one case
(see \cite[\S4]{10}, \cite[Example 2.29]{9}):
\begin{equation}
\FSE_\al(x):=\FSJ_\al(x)+\frac{i\,x}{2(\al+1)}\,\FSJ_{\al+1}(x).
\label{5}
\end{equation}
In particular,
$\FSE_{-1/2}(x)=e^{ix}$.
The nonsymmetric Hankel transform pair takes the form
\begin{equation}
\begin{cases}
&\displaystyle\wh f(\la)=\int_{-\iy}^\iy f(x)\,\FSE_\al(-\la x)\,
|x|^{2\al+1}\,dx,\mLP
&\displaystyle
f(x)=\frac1{2^{2(\al+1)}\Ga(\al+1)^2}\int_{-\iy}^\iy \wh f(\la)\,
\FSE_\al(\la x)\,|\la|^{2\al+1}\,d\la.
\end{cases}
\label{2}
\end{equation}
The transform pair \eqref{2} follows immediately from \eqref{1} by
putting $f(x)=f_1(x)+x f_2(x)$ in \eqref{2} with $f_1$ and $f_2$ even.
For given $\al$ define the differential-reflection operator
\begin{equation}
(Yf)(x):=f'(x)+(\al+\thalf)\,\frac{f(x)-f(-x)}x\,.
\label{72}
\end{equation}
This is the {\em Dunkl operator} for root system $A_1$
(see \cite[Definition 4.4.2)]{11}).
Then we have the eigenvalue equation
\begin{equation}
Y(\FSE_\al(\la\,.\,))=i\,\la\,\FSE_\al(\la\,.\,)\,.
\label{3}
\end{equation}
If, in \eqref{3}, we substitute \eqref{5}, compare even and odd parts,
and then substitute \eqref{4}, then we see that \eqref{3} is
equivalent with a pair of lowering and raising differentiation formulas
for Bessel functions (see \cite[(10.6.2)]{8}):
\[
J_\al'(x)-\,\frac\al x\,J_\al(x)=-J_{\al+1}(x),\qquad
J_{\al+1}'(x)+\frac{\al+1}x\,J_{\al+1}(x)=J_\al(x).
\]
\section{Askey-Wilson polynomials}
Askey-Wilson polynomials were introduced in \cite{12},
see also \cite[\S7.5]{1} and \cite[\S3.1]{2}.
We will consider these polynomials in $x=(z+z^{-1})/2$
as symmetric Laurent polynomials in $z$ and we will renormalize them
such that they are {\em monic}, i.e., the coefficient of the highest
degree term in $z$ is 1:
\begin{multline}
\qquad P_n[z]=P_n[z;a,b,c,d\mid q]=P_n\big(\thalf(z+z^{-1})\bigr)\\
:=
\frac{(ab,ac,ad;q)_n}{a^n(abcdq^{n-1};q)_n}
\,\qhyp43{q^{-n},q^{n-1}abcd,az,az^{-1}}{ab,ac,ad}{q,q}.\qquad
\label{13}
\end{multline}
They are eigenfunctions of a second order $q$-difference operator
$L=L_{a,b,c,d;q}$\,:
\begin{equation}
(LP_n)[z]=(q^{-n}-1)(1-abcdq^{n-1})P_n[z],
\label{9}
\end{equation}
where
\[
(Lf)[z]:=
A[z]\,f[qz]+A[z^{-1}]\,f[q^{-1}z]-\big(A[z]+A[z^{-1}]\big)\,f[z]
\]
and
\[
A[z]:=\frac{(1-az)(1-bz)(1-cz)(1-dz)}{(1-z^2)(1-qz^2)}\,.
\]
The Askey-Wilson polynomials are on the top level of the $q$-Askey scheme
(see \cite[\S3]{2}; a graphical display of this scheme including the arrows
will appear in \cite{13}, just before Chapter~14).

The three-term recurrence relation for the polynomials \eqref{13}
 is as follows (see \cite[(3.1.5)]{2}):
\begin{equation}
(z+z^{-1})P_n[z]=P_{n+1}[z]+B_n P_n[z]+C_n P_{n-1}[z],
\label{14}
\end{equation}
where $C_n P_{n-1}[z]:=0$ for $n=0$, and
where $B_n$ and $C_n$ (partially expresed in terms of the $e_i$, see
\eqref{21}) are given by:
\begin{align}
B_n&:=q^{n-1}\,\frac{(1-q^n-q^{n+1})e_3+qe_1+q^{2n-1}e_3e_4
-q^{n-1}(1+q-q^{n+1})e_1e_4}{(1-q^{2n-2}e_4)(1-q^{2n}e_4)}\,,
\label{15}\sLP
C_n&:=(1-q^{n-1}ab)(1-q^{n-1}ac)(1-q^{n-1}ad)(1-q^{n-1}bc)
(1-q^{n-1}bd)(1-q^{n-1}cd)
\nonu\\
&\qquad\qquad\times\frac{(1-q^n)(1-q^{n-2}e_4)}
{(1-q^{2n-3}e_4)(1-q^{2n-2}e_4)^2(1-q^{2n-1}e_4)}\,.
\label{16}
\end{align}
From \eqref{15} and \eqref{16} it is clear that $P_n[z;a,b,c,d\mid q]$
is symmetric in $a,b,c,d$ (well-known, but not yet evident from
\eqref{13}).

By Favard's theorem
(see \cite[Theorems I.4.4 and II.3.2]{15}), there exists a positive
Borel measure $\mu=\mu_{a,b,c,d;q}$ on $\RR$ with $\mu(\RR)=1$ such that
\begin{equation}
\lan P_m,P_n\ran=\lan P_m,P_n\ran_{a,b,c,d;q}
:=\int_\RR P_m(x)\,P_n(x)\,d\mu(x)=h_n\,\de_{m,n}
\label{17}
\end{equation}
with
\begin{equation}
h_n=h_n^{a,b,c,d;q}=C_1C_2\ldots C_n
=\frac{(q,ab,ac,ad,bc,bd,cd;q)_n}
{(abcd;q)_{2n}(q^{n-1}abcd;q)_n}>0.
\label{18}
\end{equation}
if and only if
\begin{equation}
\mbox{$B_n$ is real for $n\ge0$ and $C_n>0$ for $n\ge1$.}
\label{19}
\end{equation}
A sufficient condition for \eqref{19} to hold is that among $a,b,c,d$
there are two, one or zero pairs of complex conjugates with the other
parameters being real and that pairwise products of parameters are
less than 1 in absolute value. Then the orthogonality measure $\mu$
can be given explicitly (see \cite[\S2]{12}, \cite[\S3.1]{2}).
\section{Nonsymmetric Askey-Wilson polynomials}
\label{49}
From now on also assume:
\begin{equation}
a,b,c,d\ne0,\quad
abcd\ne q^{-m} (m=0,1,2,\ldots),\quad
\{a,b\}\cap\{a^{-1},b^{-1}\}=\emptyset.
\label{20}
\end{equation}
In terms of
\begin{align*}
P_n[z]&\;=P_n[z;a,b,c,d\mid q],\\
Q_n[z]&:=a^{-1}b^{-1}z^{-1}(1-az)(1-bz)\,P_{n-1}[z;qa,qb,c,d\mid q]
\end{align*}
we define the {\em nonsymmetric Askey-Wilson polynomials} by:
\begin{align}
E_{-n}&:=P_n-Q_n\qquad(n=1,2,\ldots),\label{25}
\\
E_n&:=P_n-\frac{ab(1-q^n)(1-q^{n-1}cd)}{(1-q^n ab)(1-q^{n-1}abcd)}\,Q_n\qquad(n=0,1,2,\ldots),
\label{26}
\end{align}
with the convention that $(1-q^n)Q_n:=0$ for $n=0$.
They are eigenfunctions of a $q$-difference-reflection operator
\begin{align}
(Yf)[z]
:=q^{-1}abcd\,f[z]
+\frac{(1-az)(1-bz)(1-cz)(1-dz)}{(1-z^2)(1-q z^2)}(f[qz]-f[z])&\nonu\sLP
+\frac{(1-a z)(1-b z) \bigl((c+d)qz-(cd+q)\bigr)}
{q(1-z^2)(1-q z^2)}\,(f[z^{-1}]-f[z])&\nonu\sLP
+\frac{(c-z)(d-z)\bigl(1+ab-(a+b)z\bigr)}{(1-z^2)(q-z^2)}\,
(f[qz^{-1}]-f[z])&.
\label{6}
\end{align}
Then
\begin{align}
YE_{-n}=q^{-n}\,E_{-n}&\qquad(n=1,2,\ldots),
\label{7}\\
YE_n=q^{n-1}abcd\,E_n&\qquad(n=0,1,2,\ldots).
\label{8}
\end{align}
The idea of working with these functions $E_n$ and operator $Y$
comes from Cherednik's theory \cite{14}
of double affine Hecke algebras associated with root systems,
extended by S.\ Sahi \cite{6}, \cite{7} to the type
$(C_l^\vee,C_l)$. The special case $l=1$ was treated afterwards in
\cite{5}, \cite[\S6.6]{4} and \cite{3}.

The {\em double affine Hecke algebra} of type $(C_1^\vee,C_1)$
is the algebra $\DAHA$ generated by $Z$, $Z^{-1}$, $T_1$, $T_0$
with relations
$ZZ^{-1}=1=Z^{-1}Z$ and
\begin{align*}
(T_1+ab)(T_1+1)&=0,&
(T_0+q^{-1}cd)(T_0+1)&=0,\\
(T_1Z+a)(T_1Z+b)&=0,&
(qT_0Z^{-1}+c)(qT_0Z^{-1}+d)&=0.
\end{align*}
This algebra acts faithfully on the space of Laurent polynomials:
\begin{align*}
(Zf)[z]&:=z\,f[z],
\sLP
(T_1f)[z]&:=\frac{(a+b)z-(1+ab)}{1-z^2}\,f[z]+
\frac{(1-az)(1-bz)}{1-z^2}\,f[z^{-1}],
\sLP
(T_0f)[z]&:=\frac{q^{-1}z((cd+q)z-(c+d)q)}{q-z^2}\,f[z]
-\frac{(c-z)(d-z)}{q-z^2}\,f[qz^{-1}].
\end{align*}
Then $Y=T_1T_0$ acts on the space of symmetric Laurent polynomials as
\eqref{6}.

The following Proposition (see \cite[Proposition 3.1]{3})
will enable us to rewrite \eqref{7} and \eqref{8}
in vector-valued form.
\begin{proposition}
\quad\\
{\rm(a)}\quad
$T_1$ acting on Laurent polynomials
has eigenvalues $-ab$ and $-1$.\\
{\rm(b)}\quad
$T_1f=-ab\,f$\quad $\,\Equiva$\quad $f$ is symmetric.\\
{\rm(c)}\quad
$T_1f=-f$ \qquad $\Equiva$\quad $f[z]=z^{-1}(1-az)(1-bz)g[z]$\\
\indent\;
for some symmetric Laurent polynomial $g$.
\end{proposition}
Let $A$ be an operator acting on the space of Laurent polynomials. Write
\begin{equation}
f[z]=f_1[z]+z^{-1}(1-az)(1-bz)f_2[z]\quad
\mbox{($f_1,f_2$ symmetric Laurent polynomials)}.
\label{24}
\end{equation}
Then we can write
\begin{equation}
(Af)[z]=
(A_{11}f_1+A_{12}f_2)[z]+z^{-1}(1-az)(1-bz)
(A_{21}f_1+A_{22}f_2)[z],
\label{53}
\end{equation}
where the $A_{i\,j}$ are operators acting on the space of symmetric
Laurent polynomials.
So we have the identifications
\begin{equation}
f\leftrightarrow\bma f_1\\f_2\ema,\quad
A\leftrightarrow\bma A_{11}&A_{12}\\A_{21}&A_{22}\ema.
\label{54}
\end{equation}
In particular, we have
\begin{align}
E_{-n}[z]&=
\bma P_n[z;a,b,c,d\mid q]\sLP-a^{-1}b^{-1} P_{n-1}[z;qa,qb,c,d\mid q]\ema
\quad(n=1,2,\ldots),
\label{27}\\
E_n[z]&=
\bma P_n[z;a,b,c,d\mid q]\sLP
\dstyle-\frac{(1-q^n)(1-q^{n-1}cd)}{(1-q^nab)(1-q^{n-1}abcd)}\,
P_{n-1}[z;qa,qb,c,d\mid q]\ema\quad
(n=0,1,2,\ldots),
\label{28}
\end{align}
where $(1-q^n)P_{n-1}:=0$ for $n=0$.
Also $Y=\bma Y_{11}&Y_{12}\\Y_{21}&Y_{22}\ema$ with
\begin{align}
Y_{11}&=q^{-1}abcd-\frac{ab}{1-ab}\,L_{a,b,c,d;q}\,,\label{33}\\
Y_{22}&=\frac{1-abcd-abq+abcdq+L_{aq,bq,c,d;q}}{q(1-ab)}\label{34}
\end{align}
(recall that $L_{a,b,c,d;q}$ is the second order $q$-difference
operator occurring in \eqref{9}) and
\begin{align}
&(Y_{21}g)[z]=
\frac{z(c-z)(d-z)\,\bigl(g[q^{-1}z]-g[z]\bigr)}
{(1-ab)(1-z^2)(q-z^2)}
+\frac{z(1-cz)(1-dz)\,\bigl(g[qz]-g[z]\bigr)}
{(1-ab)(1-z^2)(1-qz^2)}\,,\label{35}\bLP
&(Y_{12}h)[z]=
\frac{ab(a-z)(b-z)(1-az)(1-bz)}{(1-ab)z(q-z^2)(1-qz^2)}\,
\bigl((cd+q)(1+z^2)-(1+q)(c+d)z\bigr)\,h[z]\nonu\sLP
&\hskip4cm-\frac{ab(a-z)(b-z)(c-z)(d-z)(aq-z)(bq-z)}{q(1-ab)z(1-z^2)(q-z^2)}\,
h[q^{-1}z]\nonu\sLP
&\hskip4cm-\frac{ab(1-az)(1-bz)(1-cz)(1-dz)(1-aqz)(1-bqz)}
{q(1-ab)z(1-z^2)(1-qz^2)}\,h[qz].\label{36}
\end{align}

The eigenvalue equations \eqref{7}, \eqref{8} for $E_n$ and for $E_{-n}$
are equivalent to the four equations
\begin{align*}
&L_{a,b,c,d;q}P_n[\,.\,;a,b,c,d\mid q]=
(q^{-n}-1)(1-abcdq^{n-1})P_n[\,.\,;a,b,c,d\mid q],\\
&L_{qa,qb,c,d;q}P_{n-1}[\,.\,;qa,qb,c,d\mid q]=
(q^{-n+1}-1)(1-abcdq^n)P_{n-1}[\,.\,;qa,qb,c,d\mid q],\\
&Y_{21}P_n[\,.\,;a,b,c,d\mid q]
=-\,\frac{(q^{-n}-1)(1-cdq^{n-1})}{1-ab}\,
P_{n-1}[\,.\,;qa,qb,c,d\mid q],\\
&Y_{12}P_{n-1}[\,.\,;qa,qb,c,d\mid q]=
-\,\frac{ab(q^{-n}-ab)(1-abcdq^{n-1})}{1-ab}\,
P_n[\,.\,;a,b,c,d\mid q].
\end{align*}
\section{Orthogonality relations}
\label{37}
Consider $\lan\,.\,,\,.\,\ran_{a,b,c,d;q}$ (see \eqref{17}, \eqref{18})
as a symmetric bilinear form on the space of symmetric Laurent polynomials.
With the identification $f\leftrightarrow(f_1,f_2)$ between a Laurent polynomial $f$ and a pair of symmetric Laurent polynomials (see \eqref{24})
we look for a symmetric bilinear form on the space of Laurent polynomials
of the form
\begin{equation}
\lan g,h\ran=\lan(g_1,g_2),(h_1,h_2)\ran=
\lan g_1,h_1\ran_{a,b,c,d;q}+C\lan g_2,h_2\ran_{qa,qb,c,d;q}
\label{29}
\end{equation}
such that the nonsymmetric Askey-Wilson polynomials $E_n$ ($n\in\ZZ$)
given by \eqref{25}, \eqref{26} are orthogonal with respect to this form, i.e.
\begin{equation}
\lan E_m,E_n\ran=0\quad(m\ne n).
\label{30}
\end{equation}
By \eqref{27}, \eqref{28} the orthogonality
certainly holds if $|n|\ne|m|$. Thus we have to determine $C$ in \eqref{29}
such that $\lan E_n,E_{-n}\ran=0$. By \eqref{25}, \eqref{26} this turns
down to
\[
C=-ab\,\frac{(1-q^nab)(1-q^{n-1}abcd)}{(1-q^n)(1-q^{n-1}cd)}\,
\frac{h_n^{a,b,c,d;q}}{h_{n-1}^{qa,qb,c,d;q}}\,.
\]
A priori, it is not clear that $C$ is independent of $n$. However, 
from \eqref{18} we compute
\begin{equation*}
\frac{h_n^{a,b,c,d;q}}{h_{n-1}^{qa,qb,c,d;q}}=
\frac{(1-q^n)(1-q^{n-1}cd)}
{(1-q^nab)(1-q^{n-1}abcd)}\,
\frac{(1-ab)(1-qab)(1-ac)(1-ad)(1-bc)(1-bd)}
{(1-abcd)(1-qabcd)}\,.
\end{equation*}
Thus $C$ is independent of $n$ and the form \eqref{29} becomes more explicitly
\begin{equation}
\lan g,h\ran=
\lan g_1,h_1\ran_{a,b,c,d;q}
-ab\,\frac{(1-ab)(1-qab)(1-ac)(1-ad)(1-bc)(1-bd)}
{(1-abcd)(1-qabcd)}\,\lan g_2,h_2\ran_{qa,qb,c,d;q}\,.
\label{22}
\end{equation}
With respect to the form \eqref{22} we have thus shown that
the orthogonality \eqref{30} holds.

Under some further assumptions we can even show that \eqref{22} defines
a positive definite inner product.
\begin{proposition}
Let $a,b,c,d;q$, besides satisfying \eqref{20}, also be such
that condition \eqref{19} holds both for $a,b,c,d;q$ and for $qa,qb,c,d;q$.
Moreover assume that $ab<0$. Then $cd<1$, all coefficients in
\eqref{25}, \eqref{26} are real, and
the inner product defined by \eqref{22}
is positive definite.
\end{proposition}
\Proof
Since $C_n$ for $a,b,c,d;q$ and $C_{n-1}$ for $qa,qb,c,d;q$ are positive,
their quotient must be positive. Hence, by \eqref{16},
\[
\frac{(1-q^{n-1}cd)}{(1-q^{n-2}cd)}\,\frac{(1-q^{n-2}abcd)}{(1-q^{n-1}abcd)}>0
\qquad(n\ge2).
\]
Taking the product of these inequalities and telescoping yields
$\frac{1-abcd}{1-cd}>0$.
Since $ab<0$, we conclude that $cd$ is real. Then $cd\ge1$ is impossible,
so $cd<1$. Hence $abcd<1$.
Then we also see that all coefficients in
\eqref{25}, \eqref{26} are real. Finally, since
\[
C_1=\frac{(1-ab)(1-ac)(1-ad)(1-bc)(1-bd)(1-cd)(1-q)}
{(1-abcd)^2(1-qabcd)}>0,
\]
we have
\[
C=-ab\,\frac{(1-ab)(1-qab)(1-ac)(1-ad)(1-bc)(1-bd)}
{(1-abcd)(1-qabcd)}>0,
\]
and the positive definiteness of \eqref{22} is settled.\qed
\bPP
In Noumi \& Stokman \cite[Proposition 6.8]{5}
a biorthogonality result involving the system of Laurent polynomials
$E_n$ is given with respect to a bilinear form defined in terms of a
contour integral. This is closely connected with our orthogonality
\eqref{30}. Our positive definiteness result is probably new. It does not
seem to be possible to rewrite the contour integral in \cite{5}
such that a positive definite hermitian inner product will become apparent
under suitable constraints on the parameters. In the $q=1$ limit case this will
turn out to be much nicer, see \eqref{67}.
\section{From nonsymmetric Askey-Wilson to nonsymmetric little $q$-Jacobi}
We consider {\em little $q$-Jacobi polynomials} (see \cite[\S3.12]{2})
in monic form:
\begin{equation}
P_n(x;a,b;q):=\frac{(-1)^n q^{n(n-1)/2}(aq;q)_n}
{(abq^{n+1};q)_n}\,\qhyp21{q^n,abq^{n+1}}{aq}{q,qx}.
\label{31}
\end{equation}
They are limits of Askey-Wilson polynomials \eqref{13}
(see \cite[Proposition 6.3]{16}):
\begin{equation}
P_n(x;a,b;q)=\lim_{\la\downarrow0}\la^n
P_n[\la^{-1}x;-q^{1/2}a,qb\la,-q^{1/2},\la^{-1}\mid q].
\label{32}
\end{equation}
Note that we take the limit of symmetric Laurent polynomials in $x$,
but in the limit we have ordinary polynomials in $x$, since all negative
powers of $x$ are killed in the limit.

The polynomials \eqref{31} are eigenfunctions of a second order
$q$-difference operator
$L=L_{a,b;q}$\,:
\begin{equation}
(LP_n)(x)=(q^{-n}-1)(1-abq^{n+1})P_n(x),
\label{39}
\end{equation}
where
\[
(Lf)(x):=
A(x)\,f(qx)+B(x)\,f(q^{-1}x)-\big(A(x)+B(x)\big)\,f(x)
\]
and
\[
A(x):=\frac{abqx-a}x\,,\qquad
B(x):=\frac{x-1}x\,.
\]

The system of polynomials \eqref{31} is orthogonal with respect to
a positive orthogonality measure if $0<a<q^{-1}$ and $b<q^{-1}$.
This can be given as an explicit $q$-integral, see \cite[(3.12.2)]{2}.
In the orthogonality relations
\begin{equation*}
\lan P_m,P_n\ran=\lan P_m,P_n\ran_{a,b;q}=h_n\,\de_{m,n}
\end{equation*}
we have
\begin{equation*}
h_n=h_n^{a,b;q}=\frac{q^{n^2}a^n(q,aq,bq;q)_n}
{(abq^2;q)_{2n}(abq^{n+1};q)_n}\,.
\end{equation*}

Corresponding to \eqref{32}
there are limits of nonsymmetric Askey-Wilson polynomials
\eqref{27}, \eqref{28} (in vector-valued form) which yield {\em nonsymmetric
little $q$-Jacobi polynomials} in vector-valued form:
\[
E_n(x;a,b,q):=\lim_{\la\downarrow0}\la^n
E_n[\la^{-1}x;-q^{1/2}a,qb\la,-q^{1/2},\la^{-1}\mid q]\qquad(n\in\ZZ).
\]
Their expressions in vector-valued form (with the usual convention
for $n=0$) are
\begin{align}
E_{-n}(x;a,b,q)
&=\bma P_n(x;a,b;q)\sLP q^{-3/2}a^{-1}b^{-1} P_{n-1}(x;qa,qb;q)\ema
\quad(n=1,2,\ldots),
\label{42}\\
E_n(x;a,b,q)&=\bma P_n(x;a,b;q)\sLP
\dstyle-\frac{q^{n-\half}(1-q^n)}{1-q^{n+1}ab}\, P_{n-1}(x;qa,qb;q)\ema
\quad(n=0,1,2,\ldots).
\label{43}
\end{align}
Note that taking corresponding limits of the nonsymmetric Askey-Wilson
polynomials \eqref{25}, \eqref{26} (as Laurent polynomials) would
give a system of linearly dependent ordinary polynomials.

The $2\times 2$ matrix-valued operator
$Y=\bma Y_{11}&Y_{12}\\Y_{21}&Y_{22}\ema$ with entries given by
\eqref{33}--\eqref{36} also has a limit for $\la\downarrow0$ after
the rescaling $a\to-q^{1/2}a$, $b\to qb\la$, $c\to-q^{1/2}$, $d\to\la^{-1}$,
$z\to\la^{-1} x$. Thus, if the eigenvalue equations \eqref{7}, \eqref{8}
are rescaled in this way with both sides being multiplied by $\la^n$, then
in the limit for $\la\downarrow0$ we obtain
\begin{align}
&\left(\bma Y_{11}&Y_{12}\\Y_{21}&Y_{22}\ema-q^{-n}\right)
\bma P_n(x;a,b;q)\sLP q^{-3/2}a^{-1}b^{-1} P_{n-1}(x;qa,qb;q)\ema=0\qquad
(n>0),
\label{44}
\mLP
&\left(\bma Y_{11}&Y_{12}\\Y_{21}&Y_{22}\ema-q^{n+1}ab\right)
\bma P_n(x;a,b;q)\sLP
\dstyle-\,\frac{q^{n-\half}(1-q^n)}{1-q^{n+1}ab}\,
P_{n-1}(x;qa,qb;q)\ema=0\qquad(n\ge0).
\label{45}
\end{align}
Here
\begin{equation}
\label{46}
\begin{split}
&Y_{11}=qab,\qquad
Y_{22}=q^{-1}-q(1-q)ab+q^{-1}L_{aq,bq;q},\\
&(Y_{21}g)(x)=
\frac{g(x)-g(qx)}{q^{1/2}x}\,,\quad
(Y_{12}g)(x)=a^2bq^{3/2}(1-bqx)g(x)-abq^{1/2}(1-x)g(q^{-1}x).
\end{split}
\end{equation}

In a similar way as in section \ref{37} we can prove that
that the vector-valued polynomials $E_n$ are orthogonal with respect
to the symmetric bilinear form
\begin{equation}
\lan g,h\ran=\lan(g_1,g_2),(h_1,h_2)\ran=
\lan g_1,h_1\ran_{a,b;q}+
\frac{q^2a^2b(1-qa)(1-qb)}
{(1-q^2ab)(1-q^3ab)}\,\lan g_2,h_2\ran_{qa,qb;q}\,.
\label{68}
\end{equation}
This form is positive definite if $b>0$, i.e., if
$a,b\in(0,q^{-1})$.
\section{Limits to nonsymmetric Jacobi polynomials}
Consider {\em Jacobi polynomials} (see \cite[\S1.8]{2}) as
monic symmetric Laurent polynomials:
\begin{multline}
P_n[z;\al,\be]=P_n\big((z+z^{-1})/2;\al,\be\big)
:=\frac{2^{2n}(\al+1)_n}{(n+\al+\be+1)_n}\,
\hyp21{-n,n+\al+\be+1}{\al+1}{\frac{2-z-z^{-1}}4}\\
=\frac{2^{2n}n!}{(n+\al+\be+1)_n}\,P_n^{(\al,\be)}\big((z+z^{-1})/2\big).
\label{57}
\end{multline}
For $\al,\be>-1$ these polynomials satisfy the orthogonality relations
\begin{equation}
\lan P_m,P_n\ran=\lan P_m,P_n\ran_{\al,\be}=
\frac{2^{-(\al+\be+1)}\Ga(\al+\be+2)}{\Ga(\al+1)\Ga(\be+1)}
\int_{-1}^1 P_m(x)\,P_n(x)\,(1-x)^\al(1+x)^\be\,dx=
h_n\,\de_{m,n},
\label{62}
\end{equation}
where
\begin{equation*}
h_n=h_n^{\al,\be}=\frac{2^{4n}(\al+1)_n(\be+1)_n n!}
{(\al+\be+2)_{2n}(n+\al+\be+1)_n}\,.
\end{equation*}

Jacobi polynomials are limits of Askey-Wilson polynomials, or rather of
continuous $q$-Jacobi polynomials, a two-parameter subclass of
the Askey-Wilson polynomials (see \cite[(3.10.14), (5.10.2)]{2}):
\begin{equation}
P_n[z;\al,\be]=\lim_{q\uparrow1}
P_n\big[z;q^{\al+\half},-q^{\be+\half},q^{\half},
-q^{\half}\mid q\big].
\label{48}
\end{equation}
The ordering of parameters in the \RHS\ of \eqref{48}, different from
the ordering in \cite[(3.10.14)]{2} but allowed in view of the symmetry in
the parameters, was chosen in order to be able to take limits of the
formulas in section \ref{49}.
Indeed, we can now obtain nonsymmetric Jacobi polynomials as limits
\[
E_n[z;\al,\be]:=\lim_{q\uparrow1}
E_n\big[z;q^{\al+\half},-q^{\be+\half},q^{\half},
-q^{\half}\mid q\big]
\]
of
nonsymmetric Askey-Wilson polynomials \eqref{25}, \eqref{26} by
using \eqref{48}. Then
\begin{align}
E_{-n}[z;\al,\be]&=P_n[z;\al,\be]-(z-z^{-1})\,P_{n-1}[z;\al+1,\be+1]
\qquad(n=1,2,\ldots)
\label{50},\\
E_n[z;\al,\be]&=P_n[z;\al,\be]+\frac n{n+\al+\be+1}\,
(z-z^{-1})\,P_{n-1}[z;\al+1,\be+1]
\qquad(n=0,1,2,\ldots),
\label{51}
\end{align}
where $nP_{n-1}:=0$ for $n=0$.

Starting with the operator $Y=Y_{a,b,c,d;q}$ given by \eqref{6} we can also
obtain a differential reflection operator
\begin{equation*}
Y=Y_{\al,\be}:=\lim_{q\uparrow1}\;(1-q)^{-1}
\Big(Y_{q^{\al+\half},-q^{\be+\half},q^{\half},
-q^{\half};q}-1\Big)
\end{equation*}
as a limit case. Then
\begin{equation}
(Yf)[z]=-z f'[z]+\frac{\al+\be+1+(\al-\be)z}{1-z^2}\,(f[z]-f[z^{-1}])
-(\al+\be+1) f[z].
\label{52}
\end{equation}
Similarly to \eqref{7} and \eqref{8} and as a limit case of them,
the polynomials \eqref{50}, \eqref{51} are eigenfunctions of the operator
\eqref{52}:
\begin{align}
YE_{-n}=n\,E_{-n}&\qquad(n=1,2,\ldots),
\label{55}\\
YE_n=-(n+\al+\be+1)\,E_n&\qquad(n=0,1,2,\ldots).
\label{56}
\end{align}
The method of \eqref{24}, \eqref{54} to rewrite eigenvalue equations
for nonsymmetric polynomials as similar equations for vector-valued
symmetric polynomials with a matrix-valued operator, also has a limit case
here. Now the identifications in \eqref{54} are based on the rules
\begin{align}
f[z]&=f_1[z]-(z-z^{-1})f_2[z]\quad
\mbox{($f_1,f_2$ symmetric Laurent polynomials)},
\label{64}\\
(Af)[z]&=
(A_{11}f_1+A_{12}f_2)[z]-(z-z^{-1})
(A_{21}f_1+A_{22}f_2)[z].
\end{align}
Accordingly, \eqref{55}, \eqref{56} hold with
\begin{align}
E_{-n}[z]&=(P_n[z;\al,\be),P_{n-1}[z;\al+1,\be+1])\qquad(n=1,2,\ldots),
\label{58}
\\
E_n[z]&=\Big(P_n[z;\al,\be),-\frac n{n+\al+\be+1}\,
P_{n-1}[z;\al+1,\be+1]\Big)\qquad(n=0,1,2,\ldots)
\label{59}
\end{align}
and
\begin{equation}
Y=\bma-(\al+\be+1)&(z^2-1)\frac d{dz}+(\al+\be+2)(z+z^{-1})+2(\al-\be)\mLP
(1-z^{-2})^{-1}\frac d{dz}&0\ema.
\label{60}
\end{equation}
Note that, in view of \eqref{58}, \eqref{59} and \eqref{57}
the eigenvalue equations \eqref{55}, \eqref{56} are equivalent to
the well-know pair of shift operator relations
\cite[(1.8.6), (1.8.7)]{2} for Jacobi polynomials.
\begin{remark}
The operator $Y$ in \eqref{52} coincides, up to a constant term, with
Cherednik's \cite[line after (3.25)]{14} trigonometric Dunkl operator
in the case of root system $BC_1$. In the notation of
\cite[(1.2)]{17} this is the operator
\[
\frac d{dt}-\thalf(k_1+2k_2)+
\left(\frac{k_1}{1-e^{-t}}+\frac{2k_2}{1-e^{-2t}}\right)(1-s),
\]
with $s$ the reflection operator. After substitution of $z=e^{-t}$,
$k_1=\al-\be$, $k_2=\be+\thalf$ this gives $Y+\thalf(\al+\be+1)$ with
$Y$ as in \eqref{52}. Earlier than \cite{14},
Heckman \cite{18} (see also \cite[(1.12)]{17})
proposed a trigonometric generalization of the Dunkl operator
Dunkl operator which is different from Cherednik's version.
For $BC_1$ this operator becomes
\[
\frac d{dt}+\thalf
\left(k_1\frac{1+e^{-t}}{1-e^{-t}}+2k_2\frac{1+e^{-2t}}{1-e^{-2t}}\right)(1-s).
\]
After substitution of $t=2i\tha$ and division by $2i$, and with $k_1$ and
$k_2$ as before, this gives the operator considered by
Chouchene \cite[p.1]{19}.
\end{remark}

In a similar way as before we can prove
that the vector-valued polynomials $E_n$ given by \eqref{50}, \eqref{51}
are orthogonal with respect
to the symmetric bilinear form
\begin{equation}
\lan g,h\ran=\lan(g_1,g_2),(h_1,h_2)\ran=
\lan g_1,h_1\ran_{\al,\be}+
\frac{16(\al+1)(\be+1)}
{(\al+\be+2)(\al+\be+3)}\,\lan g_2,h_2\ran_{\al+1,\be+1}\,.
\label{63}
\end{equation}
Here $\lan\,.\,,\,.\,\ran_{\al,\be}$ is as in
\eqref{62}.
This form is positive definite if $\al,\be>-1$.
The inner product \eqref{63} can also be written in
terms of an integral with positive weight function.
First observe that \eqref{64} implies that
\begin{equation}
f\left[-\left(x\pm i\sqrt{1-x^2}\,\right)^2\right]=f_1(1-2x^2)\pm
4ix\sqrt{1-x^2}\,f_2(1-2x^2)\qquad(x\in[-1,1]).
\label{65}
\end{equation}
Hence, for Laurent polynomials $g,h$ we obtain from \eqref{63} and
\eqref{65} that
\begin{align}
\lan g,h\ran&=\frac{\Ga(\al+\be+2)}{\Ga(\al+1)\Ga(\be+1)}\nonu\\
&\quad\times\int_{-1}^1 g\Big[-\big(x\pm i\sqrt{1-x^2}\,\big)^2\Big]\,
\overline{h\Big[-\big(x\pm i\sqrt{1-x^2}\,\big)^2\Big]}\,
|x|^{2\al+1} (1-x^2)^\be\,dx\\
&=\frac{2^{-(2\al+2\be+4)}\Ga(\al+\be+2)}{\Ga(\al+1)\Ga(\be+1)}
\int_{-\pi}^\pi g(e^{i\tha})\,\overline{h(e^{i\tha})}\;
\Big|(1-e^{i\tha})^{\al+\half}\,(1+e^{i\tha})^{\be+\half}\Big|^2\,d\tha.
\label{66}
\end{align}
So the orthogonality of the vector-valued polynomials $E_n$ with respect
to the inner product \eqref{63} can be rewritten as the following orthogonality
for the Laurent polynomials $E_n$ given by \eqref{50}, \eqref{51}:
\begin{equation}
\int_{-1}^1 E_m\Big[-\big(x+i\sqrt{1-x^2}\,\big)^2\Big]\,
\overline{E_n\Big[-\big(x+i\sqrt{1-x^2}\,\big)^2\Big]}\,
|x|^{2\al+1} (1-x^2)^\be\,dx=0\qquad(m\ne n),
\end{equation}
or equivalently,
\begin{equation}
\int_{-\pi}^\pi E_m\big[e^{i\tha}\,\big]\;
\overline{E_n\big[e^{i\tha}\,\big]}\;\;
\Big|(1-e^{i\tha})^{\al+\half}\,(1+e^{i\tha})^{\be+\half}\Big|^2\,d\tha=0
\qquad(m\ne n).
\label{67}
\end{equation}
In the form \eqref{67} the orthogonality relations for
nonsymmetric Jacobi polynomials are the specialization to root system
$BC_1$ of Opdam's orthogonality relations \cite[Definition 2.1]{21} for
nonsymmetric Jacobi polynomials associated with general root systems.
The nonsymmetric Jacobi polynomials which occur as eigenfunctions of
Heckman's trigonometric Dunkl operator in case $BC_1$ are also
orthogonal with respect to the inner product \eqref{66},
see Chouchene \cite[(2.41)]{19}.

Jacobi polynomials can also be obtained as limits of
little $q$-Jacobi polynomials \eqref{31}
(see \cite[(5.12.1)]{2}). Define
\begin{equation}
\wt P_n(x;\al,\be):=\lim_{q\uparrow1}P_n(x;q^\al,q^\be;q).
\label{41}
\end{equation}
These are
monic Jacobi polynomials with orthogonality interval rescaled to $[0,1]$:
\begin{equation}
\wt P_n(x;\al,\be)=\frac{(-1)^n(\al+1)_n}{(n+\al+\be+1)_n}\,
\hyp21{-n,n+\al+\be+1}{\al+1}x.
\label{38}
\end{equation}
By comparing \eqref{38} with \eqref{57} we see that
\begin{equation}
P_n[z;\al,\be]=(-1)^n\,2^{2n}\,\wt P_n\Big(\frac{2-z-z^{-1}}4;\al,\be\Big).
\label{61}
\end{equation}

Corresponding to \eqref{41} we can take limits of formulas
\eqref{42}--\eqref{68}
involving nonsymmetric little $q$-Jacobi polynomials in vector-valued
form. Thus we will arrive at formulas which are essentially the same
as \eqref{58}--\eqref{63} and which can be fully identified with each other
by using \eqref{61}. However the expressions \eqref{50}, \eqref{51}
and \eqref{52} involving Laurent polynomials cannot be obtained as limits
of formulas for little $q$-Jacobi polynomials, since such formulas are
missing there.
\section{From nonsymmetric Jacobi to nonsymmetric Bessel}
In this final section we come back to the nonsymmetric Bessel functions,
discussed in section \ref{69},
and we show that they are limits of nonsymmetric Jacobi polynomials.
There is a well-known limit from Jacobi polynomials to Bessel functions,
see \cite[(18.11.5)]{8}, which we can rewrite in terms of the notations
\eqref{57} and \eqref{4} as
\begin{equation}
\lim_{n\to\iy}\frac{2^{\al+\be}\Ga(\al+1)}{\pi^{\half}n^{\al+\half}}\,
P_n\left(1-\frac{\la^2x^2}{2n^2};\al,\be\right)=\FSJ_\al(\la x).
\label{70}
\end{equation}
Now from \eqref{50}, \eqref{51}, rewritten in the form of \eqref{65},
we see that \eqref{70} implies:
\begin{equation}
\lim_{n\to\iy}\frac{2^{\al+\be}\Ga(\al+1)}{\pi^{\half}n^{\al+\half}}
E_{\pm n}\Big[-\Big(\tfrac{\la x}{2n}\pm i
\sqrt{1-\tfrac{\la^2 x^2}{4n^2}}\,\Big)^2\,\Big]=
\FSJ_\al(\la x)\pm\frac{i\la x}{2(\al+1)}\,\FSJ_{\al+1}(\la x)
=\FSE_\al(\pm\la x).
\label{71}
\end{equation}
Corresponding with \eqref{71} there are limits from \eqref{52}, \eqref{55},
\eqref{56} to \eqref{72} and \eqref{3}.

Finally, the limit case of \eqref{58}, \eqref{59}, \eqref{60}
turns down to \eqref{3} rewritten in vector-valued form as
\begin{equation}
\bma
0&x\frac{d}{dx}+2(\al+1)\\\frac{1}{x}\frac{d}{dx}&0\ema
\bma\mathcal{J_{\al}}(\lambda x)\\
\frac{i\lambda}{2(\al+1)}\mathcal{J}_{\al +1}(\lambda x)\ema=i\la
\bma\mathcal{J_{\al}}(\lambda x)\\
\frac{i\lambda}{2(\al+1)}\mathcal{J}_{\al +1}(\lambda x)\ema.
\end{equation}

\end{document}